\documentclass[journal]{IEEEtran}
\usepackage{amsmath}    
\usepackage{graphics,epsfig,subfigure}    
\usepackage{verbatim}   
\usepackage{color}      
\usepackage{subfigure}  
\usepackage{hyperref}   
\usepackage{amssymb}
\usepackage{epstopdf}
\usepackage{graphicx}
\usepackage{mathrsfs} 

\newcommand{\prob}{Pr}
\newcommand{\expect}[1]{\mathbb{E}\big\{#1\big\}}

\newcommand{\bv}[1]{{\boldsymbol{#1} }}
\newcommand{\script}[1]{{{\cal{#1} }}}

\newcommand{\rmax}{\text{argmax}}
\begin{document}

\title{Utility Optimal Scheduling in Energy Harvesting Networks}
 
\author{\large{Longbo Huang,  Michael J. Neely }  
\thanks{Longbo Huang,    (web:  http://www-scf.usc.edu/$\sim$longbohu) 
and Michael J. Neely (web:  http://www-rcf.usc.edu/$\sim$mjneely)
are with the Dept. of Electrical
Engineering, University of Southern California, Los Angeles, CA 90089, USA.} 
\thanks{This material is supported in part  by one or more of 
the following: the DARPA IT-MANET program
grant W911NF-07-0028, the NSF grant OCE 0520324, 
the NSF Career grant CCF-0747525.} }

\maketitle
\thispagestyle{empty}
\pagestyle{empty}
\newtheorem{remark}{Remark}
\newtheorem{fact_def}{\textbf{Fact}}
\newtheorem{coro}{\textbf{Corollary}}
\newtheorem{lemma}{\textbf{Lemma}}
\newtheorem{main}{\textbf{Proposition}}
\newtheorem{theorem}{\textbf{Theorem}}
\newtheorem{claim}{\emph{Claim}}
\newtheorem{prop}{Proposition}
\newtheorem{property}{Property}
\newtheorem{assumption}{\textbf{Assumption}}

\begin{abstract}
In this paper, we show how to achieve close-to-optimal utility performance in energy harvesting networks with only finite capacity energy storage devices. 
In these networks, nodes are capable of harvesting energy from the environment. 
The amount of energy that can be harvested  is  time varying and evolves according to some probability law. 
We develop an \emph{online} algorithm, called the Energy-limited Scheduling Algorithm (ESA),  which jointly manages the energy and makes power allocation decisions for packet transmissions. ESA only has to keep track of the amount of energy left at the network nodes and \emph{does not require any knowledge} of the harvestable energy process. We show that ESA achieves a utility that is within $O(\epsilon)$ of the optimal, for any $\epsilon>0$, while ensuring that the network congestion and the required capacity of the energy storage devices are \emph{deterministically} upper bounded by bounds of size $O(1/\epsilon)$. We then also develop the Modified-ESA algorithm (MESA) to achieve the same $O(\epsilon)$ close-to-utility performance, with the average network congestion and the required capacity of the energy storage devices being only $O([\log(1/\epsilon)]^2)$. 
\end{abstract}

\section{Introduction}
Recent developments in hardware design have enabled many general wireless networks to support themselves by harvesting energy from the environment. 
For instance, by converting mechanical  vibration into energy \cite{vibrate-to-energy-01}, by using solar panels \cite{solar-energy-05}, by utilizing  thermoeletric generators \cite{energy-harvesting-survey08}, or by converting ambient radio power into energy \cite{enhant-mobicom09}. Such harvesting methods are also referred to as  ``recycling'' energy \cite{energy-harvesting-economist10}. 
This energy harvesting ability is crucial for many network design problems. It  frees the network devices from having an ``always on'' energy source and provides a way of operating the network with a potentially infinite lifetime. 
These two advantages are particularly useful for networks that work autonomously, e.g., wireless sensor networks that perform monitoring tasks in dangerous fields \cite{sensor-volcano06}, tactical networks \cite{barrage-relay-ita10},  or wireless handheld devices that operate over a longer period \cite{energy-harvesting-mitreview09}, etc. 

However, to take full advantage of the energy harvesting technology,   efficient scheduling algorithms must consider the finite capacity for energy storage at each network node. 
In this paper, we consider the problem of constructing utility optimal scheduling algorithms in a discrete stochastic network, where the communication links have time-varying qualities, and the nodes are powered by finite capacity energy storage devices but are capable of harvesting energy. 
Every time slot, the network decides how much new data to admit and how much power to allocate over each communication link for data transmission. 
The objective of the network is to  maximize the aggregate traffic utility subject to the  constraint that the average network backlog is finite, and the ``energy-availability'' constraint is met, i.e., at all time, the energy consumed is no more than the energy stored.  
We see that the  ``energy-availability'' constraint greatly complicates the design of an efficient scheduling algorithm, due to the fact that the current energy expenditure decision may cause energy outage in the future and thus affect the future decisions. Such problems can in principle be formulated as dynamic programs (DP) and be solved optimally. However, the DP approach typically requires substantial statistical knowledge of the harvestable energy process and the channel state process, and often runs into the ``curse-of-dimensionality'' problem when the network size is large.

There has been many previous works developing algorithms for such energy harvesting networks. \cite{opt-energy-twc10} develops algorithms for a single sensor node for achieving maximum capacity and minimizing delay when the rate-power curve is linear. \cite{power-harvesting-kansal07} considers the problem of optimal power management for sensor nodes, under the assumption that the harvested energy satisfies a leaky-bucket type property. \cite{energy-tradeoff-10} looks at the problem of designing energy-efficient schemes for maximizing the decay exponent of the queue length. \cite{control-rechargeable-twc10} develops scheduling algorithms to achieve close-to-optimal utility for energy harvesting networks with time varying channels.  
\cite{power-routing-lin-infocom05} develops an energy-aware routing scheme that  approaches optimal as the network size increases.  
Outside the energy harvesting context, \cite{power-routing-lin-adhocnet07} considers the problem of maximizing the lifetime of a network with finite energy capacity and constructs a scheme that achieves a close-to-maximum lifetime. 
\cite{neelyenergy} and \cite{neelyenergydelay} develop algorithms for minimizing the time average network energy consumption for stochastic networks with ``always on'' energy source. However,  most of the existing results focus on single-hop networks and often require sufficient statistical knowledge of the harvestable energy, and results for multihop networks often do not give explicit queueing bounds and do not provide explicitly characterizations of  the needed energy storage capacities. 

We  tackle this problem using the Lyapunov optimization technique developed in \cite{neelyenergy} and \cite{neelynowbook}, combined with the idea of weight perturbation, e.g., \cite{neelyhuang_assembly} and \cite{huangneely-spn}. 
The idea of this approach is to construct the algorithm based on a quadratic Lyapunov function, but carefully perturb the weights used for decision making, so as to ``push'' the target queue levels towards certain nonzero values to avoid underflow (in our case, the target queue levels are the energy levels at the nodes). 
Based on this idea, we construct our Energy-limited Scheduling Algorithm (ESA) for achieving optimal utility in general multihop energy harvesting networks powered by finite capacity energy storage devices. ESA is an \emph{online} algorithm which makes greedy decisions every time slot \emph{without requiring any knowledge of the harvestable energy} and \emph{without requiring any statistical knowledge of the channel qualities}. We show that the ESA algorithm is able to achieve an average utility that is within  $O(\epsilon)$ of the optimal for any $\epsilon>0$, and only requires energy storage devices that are of $O(1/\epsilon)$ sizes. We also explicitly compute the required storage capacity and show that ESA also guarantees that the network backlog is deterministically bounded by $O(1/\epsilon)$. 
Furthermore, we develop the Modified-ESA  algorithm (MESA) to achieve the same $O(\epsilon)$ close-to-optimal utility performance with energy storage devices that are  only of $O([\log(1/\epsilon)]^2)$ sizes. 
We note that the approach of using perturbation in Lyapunov algorithms is novel. It not only allows us to resolve the energy outage problem easily, but also enables an easy analysis of the algorithm performance. 


Our paper is mostly related to the recent work \cite{control-rechargeable-twc10}, which considers a similar problem.  \cite{control-rechargeable-twc10}  uses a similar Lyapunov optimization approach (without perturbation) for algorithm design, and 
 achieves a similar $[O(\epsilon), O(1/\epsilon)]$ utility-backlog performance using energy storage sizes of $O(1/\epsilon)$ for single-hop networks. Multihop networks are also considered in \cite{control-rechargeable-twc10}. However, the performance bounds for multihop networks are given in terms of unknown parameters. 
In our paper, we compute the explicit $O(1/\epsilon)$ capacity requirements for the data buffers and energy storage devices for general multihop networks for achieving the $O(\epsilon)$ close-to-optimal utility performance. 
We then also  develop a scheme to achieve the same utility performance  with only $O([\log(1/\epsilon)]^2)$ energy storage capacities. 

Our paper is organized as follows: In Section \ref{section:model} we state our network model and the objective. In Section \ref{section:upperbound} we first derive an upper bound on the maximum utility. Section \ref{section:alg} presents the ESA algorithm. The $[O(\epsilon), O(1/\epsilon)]$ performance results of the ESA algorithm are presented in Section \ref{section:analysis}, for both the cases when the network randomness is i.i.d. and Markovian. We then construct the Modified-ESA algorithm (MESA) to achieve the same $O(\epsilon)$ close-to-optimal utility performance with only $O([\log(1/\epsilon)]^2)$ energy storage sizes in Section \ref{section:reduce-buffer}. 
Simulation results are presented in Section \ref{section:simulation}. We conclude our paper in Section \ref{section:conclusion}.

\section{The Network Model}\label{section:model}
We consider a general interconnected network that operates in slotted time. The network is modeled by a directed graph $\script{G}=(\script{N}, \script{L})$, where $\script{N}=\{1, 2, ..., N\}$ is the set of the $N$ nodes in the network, and $\script{L}=\{[n, m], \,\, n, m\in\script{N}\}$ is the set of communication links in the network. For each node $n$, we use $\script{N}_n^{(o)}$ to denote the set of nodes $b$ with $[n, b]\in\script{L}$, and use $\script{N}_n^{(in)}$ to denote the set of nodes $a$ with $[a, n]\in\script{L}$. We then define $d_{max}\triangleq\max_n|\script{N}^{(in)}_n|$
to be the maximum in-degree that any node $n\in\script{N}$ can have. 



\subsection{The Traffic and Utility Model}

At every time slot, the network decides how many packets destined for node $c$ to admit at node $n$. We call these traffic the \emph{commodity $c$} data and use $R^{(c)}_n(t)$ to denote the amount of new commodity $c$ data admitted. We assume that $0\leq R^{(c)}_n(t)\leq R_{max}$ for all $n, c$ with some finite $R_{max}$ at all time. \footnote{Note that this setting implicitly assumes that nodes always have packets to admit. The case when the number of packets available is random can also be incorporated into our model and solved by introducing auxiliary variables, as in \cite{neelyfairness}.}
We assume that each commodity is associated with a utility function $U^{(c)}_{n}(\overline{r}^{nc})$, where $\overline{r}^{nc}$ is the time average rate of the commodity $c$ traffic admitted into node $n$, defined as $\overline{r}^{nc}=\liminf_{t\rightarrow\infty}\frac{1}{t}\sum_{\tau=0}^{t-1}\expect{R^{(c)}_n(\tau)}$.  
Each $U^{(c)}_{n}(r)$ function is assumed to be increasing, continuously differentiable, and strictly concave in $r$ with a bounded first derivative and $U^{(c)}_{n}(0)=0$.  We  use $\beta^{nc}$ to denote the maximum first derivative of $U^{(c)}_{n}(r)$, i.e., $\beta^{nc}=(U^{(c)}_{n})'(0)$ and denote $\beta=\max_{n, c}\beta^{nc}$. 
 

\subsection{The Transmission Model}
In order to deliver the data to their destinations, each node needs to allocate power to each link for data transmission at every time slot. To model the effect that the transmission rates typically also depend on the link conditions and that the link conditions may be time varying, we let $\bv{S}(t)$ be the network \emph{channel state}, i.e.,  the $N$-by-$N$ matrix where the $(n, m)$ component of $\bv{S}(t)$ denotes the channel condition between nodes $n$ and $m$. We assume that $\bv{S}(t)$ takes values in some finite set $\script{S}=(s_1, ..., s_M)$. In the following, we first assume that $\bv{S}(t)$ is i.i.d. every time slot and use $\pi_{s_i}$ to denote $\prob(\bv{S}(t)=s_i)$. We will later extend our results to the case when $\bv{S}(t)$ is Markovian. 
At every time slot, if $\bv{S}(t)=s_i$, then the power allocation vector $\bv{P}(t)=(P_{[n, m]}(t), [n, m] \in\script{L})$, where $P_{[n, m]}(t)$ is the power allocated to link $[n, m]$ at time $t$, must be chosen from some feasible power allocation set $\script{P}^{(s_i)}$. We assume that $\script{P}^{(s_i)}$ is compact for all $s_i$, and that every power vector in $\script{P}^{(s_i)}$ satisfies the constraint that for each node $n$, $0\leq \sum_{b\in\script{N}^{(o)}_n}P_{[n, b]}(t)\leq P_{max}$ for some $P_{max}<\infty$. Also, we assume that setting any $P_{[n, m]}$ in a vector $\bv{P}\in\script{P}^{(s_i)}$ to zero yields another power vector that is still in $\script{P}^{(s_i)}$. 
Given the channel state $\bv{S}(t)$ and the power allocation vector $\bv{P}(t)$, the transmission rate over the link $[n, m]$ is given by the rate-power function $\mu_{[n, m]}(t)=\mu_{[n, m]}(\bv{S}(t),\bv{P}(t))$. For each $s_i$, we assume that the function $\mu_{[n, m]}(s_i, \bv{P})$ satisfies the following 
 properties:
\begin{property} \label{property1}
For any vectors $\bv{P}, \bv{P}'\in\script{P}^{(s_i)}$, where $\bv{P}'$ is obtained by changing any single component $P_{[n, m]}$ in $\bv{P}$ to zero, we have: 
\begin{eqnarray}
\mu_{[n, m]}(s_i, \bv{P})\leq \mu_{[n, m]}(s_i, \bv{P}') + \delta P_{[n, m]},\label{eq:rate-property1}
\end{eqnarray}
for some finite constant $\delta>0$. 
\end{property}
\begin{property} \label{property2}
If $\bv{P}'$ is obtained by setting the entry $P_{[n, b]}$ in $\bv{P}$ to zero, then: 
\begin{eqnarray}
\mu_{[a, m]}(s_i, \bv{P})\leq \mu_{[a, m]}(s_i, \bv{P}'),\,\,\forall\,[a, m]\neq[n, b].\label{eq:rate-property2}
\end{eqnarray}
\end{property}
Property \ref{property1} states that the rate obtained over a link $[n, m]$ is upper bounded by some linear function of the power allocated to it, whereas Property \ref{property2} states that reducing the power over any link does not reduce the rate over any other links. 
We see that Property \ref{property1} and \ref{property2} can usually be satisfied by most rate-power functions, e.g.,  when the rate function is differentiable and has finite directional derivatives with respect to power \cite{neelyenergy}, and the links do not interfere with each other. 

We also assume that there exists some finite constant $\mu_{max}$ such that $\mu_{[n, m]}(t)\leq\mu_{max}$ for all time under any power allocation vector and any  channel state $\bv{S}(t)$. 
\footnote{Note that in our transmission model, we did not explicitly take into account the reception power. However, this can easily be incorporated into our model at the expense of more complicated notations. All the results in this paper will still hold in this case.}
In the following, we also use $\mu_{[n, b]}^{(c)}(t)$ to denote the rate allocated to the commodity $c$ data over link $[n, b]$ at time $t$. It is easy to see that at any time $t$, we have: 
\begin{eqnarray}
\sum_{c}\mu_{[n, b]}^{(c)}(t)\leq\mu_{[n, b]}(t), \quad\forall\,\, [n, b]. \label{eq:ratecond}
\end{eqnarray}


\subsection{The Energy Queue Model}
We now specify the energy model. Every node in the network is assumed to be powered by a \emph{finite} capacity energy storage device, e.g., a battery or ultra-capacitor \cite{opt-energy-twc10}. We model such a device using an  \emph{energy queue}. We use the energy queue size at node $n$ at time $t$, denoted by $E_n(t)$, to measure the amount of the energy left in the storage device at node $n$ at time $t$. 
We assume that at every time, the nodes are capable of tracking its current energy level $E_n(t)$. 
In any time slot $t$, the power allocation vector $\bv{P}(t)$ must satisfy the following ``energy-availability'' constraint: 
\begin{eqnarray}
\sum_{b\in\script{N}^{(o)}_n}P_{[n, b]}(t) \leq E_n(t),\quad\forall\,\, n. \label{eq:energycond}
\end{eqnarray}
That is, the consumed power must be no more than what is available. 
Each node in the network is assumed to be capable of harvesting energy from the environment, using, for instance,  solar panels \cite{opt-energy-twc10}. 
However, the amount of harvestable energy in a time slot is typically not fixed and varies over time. We use $h_n(t)$ to denote the amount of harvestable energy by node $n$ at time $t$, and  denote $\bv{h}(t)=(h_1(t), ..., h_{N}(t))$ the harvestable energy vector at time $t$, called the \emph{energy state}.  
We assume that $\bv{h}(t)$ takes values in some finite set $\script{H}=\{\bv{h}_1, ..., \bv{h}_K\}$, and that $\bv{h}(t)$ is i.i.d. over each slot. However, components in each $\bv{h}_i$ vector may be correlated. We will later consider the case when $\bv{h}(t)$ is Markovian. In both cases, we assume that $\bv{h}(t)$ is independent of $\bv{S}(t)$. \footnote{This is for the ease of presentation. The results in this paper still hold if they are correlated.}

We let $\pi_{\bv{h}_i}=\prob(\bv{h}(t)=\bv{h}_i)$. 
We assume that there exists $h_{max}<\infty$ such that $h_n(t)\leq h_{max}$ for all $n, t$, and 
the energy harvested at time $t$ is assumed to be available for use in time $t+1$. 
In the following, it is convenient for us to assume that each energy queue has infinite capacity, and that each node can decide whether or not to harvest energy on each slot. 
We model this harvesting decision by using $e_n(t) \in[0, h_n(t)]$ to denote the amount of energy that is actually harvested at time $t$. We will show later that our algorithm always harvests energy when the energy queue is below a finite threshold of size $O(1/\epsilon)$ and drops it otherwise, thus can be implemented with finite capacity storage devices.

\subsection{Queueing Dynamics}
Let $\bv{Q}(t)=(Q^{(c)}_n(t), n, c\in\script{N})$, $t=0, 1,2, ...$ be the data queue backlog vector in the network, where $Q^{(c)}_n(t)$ is  the amount of commodity $c$ data queued at node $n$. We assume the following queueing dynamics:
\begin{eqnarray}
Q^{(c)}_n(t+1) &\leq& \big[Q^{(c)}_n(t) - \sum_{b\in\script{N}^{(o)}_n}\mu_{[n, b]}^{(c)}(t)\big]^+ \label{eq:Qdynamic}\\
&&\qquad\qquad\,\,\,+ \sum_{a\in\script{N}^{(in)}_n}\mu_{[a, n]}^{(c)}(t)+R^{(c)}_n(t),\nonumber
\end{eqnarray}
with $Q^{(c)}_n(0)=0$ for all $n, c\in\script{N}$, $Q^{(c)}_c(t)=0$ $\forall\,t$, and $[x]^+=\max[x, 0]$. The inequality in (\ref{eq:Qdynamic}) is due to the fact that some nodes may not have enough commodity $c$ packets to fill the allocated rates. 
In this paper, we say that the network is \emph{stable} if the following condition is met:
\begin{eqnarray}
\overline{\bv{Q}} \triangleq \limsup_{t\rightarrow\infty}\frac{1}{t}\sum_{\tau=0}^{t-1}\sum_{n,c} \expect{Q_n^{(c)}(\tau)} <\infty.\label{eq:queuestable}
\end{eqnarray}

Similarly, let $\bv{E}(t)=(E_n(t), n\in\script{N})$ be the vector of the energy queue sizes.  Due to the energy availability constraint (\ref{eq:energycond}), we see that for each node $n$, the energy queue $E_n(t)$ evolves  according to the following: \footnote{Note that we do not explicitly consider energy leakage due to the imperfectness of the energy storage devices. This is a valid assumption if the rate of energy leakage is very small compared to the amount spent in each time slot. }
\begin{eqnarray}
\hspace{-.3in}&&E_n(t+1) =  E_n(t) - \sum_{b\in\script{N}^{(o)}_n}P_{[n, b]}(t) +e_n(t), \label{eq:Edynamic} 
\end{eqnarray}
with $E_n(0)=0$ for all $n$. \footnote{We can also pre-store energy in the energy queue and initialize $E_n(0)$ to any finite positive value up to its capacity. The results in the paper will not be affected. }
Note again that by using the queueing dynamic (\ref{eq:Edynamic}), we start by  assuming that each energy queue has infinite capacity. Later we will show  that under our  algorithms, all the $E_n(t)$ values are \emph{determinstically} upper bounded, thus 
we only need a finite energy capacity in algorithm implementation. 

\subsection{Utility Maximization with Energy Management}
The goal of the network is thus to design a joint flow control, routing and scheduling, and energy  management algorithm that at every time slot,  admits  the right amount of data $R_n^{(c)}(t)$, chooses power allocation vector $\bv{P}(t)$ subject to (\ref{eq:energycond}), and transmits packets accordingly, so as to maximize the utility function: 
\begin{eqnarray}
U_{tot}(\overline{\bv{r}}) &=& \sum_{n, c} U^{(c)}_n(\overline{r}^{nc}),
\end{eqnarray}
subject to the network stability constraint (\ref{eq:queuestable}). Here $\overline{\bv{r}}=(\overline{r}^{nc}, \forall\, n, c\in\script{N})$ is the vector of the average expected admitted rates. Below, we will refer to this problem as the \emph{Utility Maximization with Energy Management problem} (UMEM).

\subsection{Discussion  of the Model}
Our model is quite general and can be used to model many networks where nodes are powered by finite capacity batteries. For instance, a field monitoring sensor network  \cite{sensor-volcano06}, or many mobile ad hoc networks \cite{mobile-adhoc-icm06}. 
Also, our model allows the harvestable energy to be correlated among network nodes. This is particularly useful, as in practice, nodes that are collocated 
may have similar  harvestable energy conditions. 


The main difficulty in designing an optimal scheduling policy here is imposed by the constraint (\ref{eq:energycond}). Indeed, (\ref{eq:energycond}) couples the current power allocation action and the future actions, in that a current action may cause the energy queue to be empty and hence block some power allocation actions  in the future. Problems of this kind usually have to be modeled as dynamic programs \cite{bertsekasdptbook}. However, this approach typically requires significant statistical knowledge of the network randomness, including the channel state and the energy state.  
Another way to utilize the harvested energy efficiently is by developing efficient sleep-wake policies, e.g., \cite{sleep-wake-icc09}. Although our model does not consider this aspect, our algorithm can also be used together with given sleep-wake policies to achieve good utility performance in that context. 



\section{Upper bounding the optimal network utility}\label{section:upperbound}
In this section, we first obtain an upper bound on the optimal utility. This upper bound will be useful for our later analysis. The result is presented in the following theorem, in which we use $\bv{r}^*$ to denote the optimal solution of the UMEM problem, subject to the constraint that the network nodes are powered by finite capacity energy storage devices. The $V$ parameter in the theorem can be any positive constant that is greater or equal to $1$, and is included for our later analysis. 
\begin{theorem}\label{theorem:upper-bound}
The optimal network utility $U_{tot}(\bv{r}^*)$  satisfies the following: 
\begin{eqnarray}
VU_{tot}(\bv{r}^*) \leq\phi^*,
\end{eqnarray}
where $\phi^*$ is the optimal value of the following optimization problem:
\begin{eqnarray}
\hspace{-.3in}&& \max: \, \phi = V\sum_{n, c}\sum_{k=1}^{K} \vartheta_kU^{(c)}_n(r^{nc}_k)\label{eq:opt-obj}\\
\hspace{-.3in}&& s.t.  \quad\,  \sum_{k=1}^{K} \vartheta_k r^{nc}_{k}+\sum_{s_i}\pi_{s_i} \sum_{k=1}^{K}\varrho^{(s_i)}_k\sum_{a\in\script{N}^{(in)}_n}\mu^{(c)}_{[a, n]}(s_i, \bv{P}_k^{(s_i)}) \nonumber\\
\hspace{-.3in}&&\quad \leq \sum_{s_i}\pi_{s_i} \sum_{k=1}^{K}\varrho^{(s_i)}_k\sum_{b\in\script{N}^{(o)}_n}\mu^{(c)}_{[n, b]}(s_i, \bv{P}_k^{(s_i)}), \forall\, (n,c), \label{eq:opt-rate}
\end{eqnarray}
\begin{eqnarray}
\hspace{-.3in}&& \qquad\quad\,\sum_{s_i}\pi_{s_i}\sum_{k=1}^{K} \varrho^{(s_i)}_k\sum_{b\in\script{N}^{(o)}_n} P^{(s_i)}_{k, [n, b]}  \label{eq:opt-energy}\\
\hspace{-.3in}&&  \qquad\qquad\qquad\qquad\qquad =  \sum_{\bv{h}_i} \pi_{\bv{h}_i} \sum_{k=1}^{K}\varphi^{(\bv{h}_i)}_k e_{n, k}^{(\bv{h}_i)},\forall n,\nonumber\\
\hspace{-.3in}&& \qquad \quad \bv{P}_k^{(s_i)}\in\script{P}^{(s_i)}, 0\leq\vartheta_k^{(s_i)}, \varrho^{(s_i)}_k, \varphi^{(\bv{h}_i)}_k\leq1, \forall s_i, k, \bv{h}_i, \nonumber\\
\hspace{-.3in}&& \qquad\quad \sum_{k=1}^{K}\vartheta_k=1, \sum_{k=1}^{K}\varrho^{(s_i)}_k=1, \sum_{k=1}^{K}\varphi^{(\bv{h}_i)}_k=1, \forall s_i, \bv{h}_i, \nonumber\\
\hspace{-.3in}&& \qquad \quad 0\leq r^{nc}_{k} \leq R_{max},\forall\, (n, c), \nonumber \\
\hspace{-.3in}&& \qquad \quad0\leq e_{n, k}^{(\bv{h}_i)} \leq h_{n}^{(\bv{h}_i)},\,\,\forall\,\, n, k, \bv{h}_i.\nonumber
\end{eqnarray}
Here $K=N^2+N+2$. \footnote{The number $K$ is due to the use of Caratheodory's Theorem in the proof argument used in \cite{huangneely-spn}. } $\{r^{nc}_{k}\}_{k=1}^{K}$ denotes the set of admission decisions used for each commodity flow. $\{\bv{P}_k^{(s_i)}\}_{k=1}^{K}$ denotes the set of power allocation vectors that are used when $\bv{S}(t)=s_i$.  $\mu^{(c)}_{[b, n]}(s_i, \bv{P}^{(s_i)}_k)$ is the rate allocated to commodity $c$ over link $[b, n]$ under $s_i$ and $\bv{P}^{(s_i)}_k$. $P^{(s_i)}_{k, [n, m]}$ is the power allocated to link $[n, m]$ under $\bv{P}^{(s_i)}_k$.  $\{e_{n, k}^{(\bv{h}_i)}\}_{k=1}^K$ is the set of  energy harvesting decisions of node $n$ when the energy state is $\bv{h}_i$, and $h^{(\bv{h}_i)}_n$ is the amount of harvestable energy for node $n$ when $\bv{h}(t)=\bv{h}_i$. 
\end{theorem}
\begin{proof}
The proof argument is similar to the one used in \cite{huangneely-spn}, hence is omitted for brevity. 
\end{proof}

Note that Theorem \ref{theorem:upper-bound} indeed holds under more general ergodic $S(t)$ and $\bv{h}(t)$ processes, e.g., when $S(t)$ and $\bv{h}(t)$ evolve according to some finite state irreducible and aperiodic  Markov chains. 
Also note that the objective function is not of the same form as $U_{tot}(\cdot)$. However, it can be shown, using Jensen's inequality, that the optimal value of the above optimization problem remains the same if we push $\sum_k\vartheta_k$ inside the function $U^{(c)}_n$, i.e., change the objective to 
$V\sum_{n, c}U^{(c)}_n(\sum_{k=1}^{K} \vartheta_kr^{nc}_k)$. 
Below, we first have the following lemma regarding the dual problem of (\ref{eq:opt-obj}): 
\begin{lemma}\label{lemma:duality}
The dual problem of (\ref{eq:opt-obj}) is given by:
\begin{eqnarray}
\min: \,\,\, g(\bv{\upsilon}, \bv{\nu}), \quad s.t.\,\,\, \bv{\upsilon}\succeq\bv{0},  \bv{\nu}\in\mathbb{R}^N,\label{eq:dual}
\end{eqnarray}
where $\bv{\upsilon}=(\upsilon^{(c)}_n, \forall\, (n, c))$, $\bv{\nu}=(\nu_n, \forall\,n)$ and 
$g(\bv{\upsilon}, \bv{\nu})$ is the dual function defined:
\begin{eqnarray}
\hspace{-.35in}&&g(\bv{\upsilon}, \bv{\nu}) = \sup_{r^{nc}, \bv{P}^{(s_i)}, e^{(\bv{h}_j)}_{n}}\sum_{s_i}\pi_{s_i}\sum_{\bv{h}_j}\pi_{\bv{h}_j}\bigg\{ V\sum_{n, c} U^{(c)}_n(r^{nc})  \nonumber\\
\hspace{-.35in}&&\qquad \qquad \quad   - \sum_{n}\upsilon^{(c)}_n\big[r^{nc} +\sum_{a\in\script{N}^{(in)}_n} \mu^{(c)}_{[a, n]}(s_i, \bv{P}^{(s_i)})\label{eq:dual-function} \\
\hspace{-.3in}&&\qquad \qquad \qquad\qquad\qquad \qquad    - \sum_{b\in\script{N}^{(o)}_n}\mu^{(c)}_{[n, b]}(s_i, \bv{P}^{(s_i)})\big] \nonumber\\
\hspace{-.35in}&& \qquad \qquad \qquad \qquad\qquad  - \sum_{n}\nu_n\big[ \sum_{b\in\script{N}^{(o)}_n} P_{[n, b]}^{(s_i)} -e_n^{(\bv{h}_j)} \big] \nonumber
\bigg\}.
\end{eqnarray}
Moreover, let $(\bv{\upsilon}^*, \bv{\nu}^*)$ be an optimal solution of (\ref{eq:dual}), then  
$\phi^*\leq g(\bv{\upsilon}^*, \bv{\nu}^*)$. 
\end{lemma}
\begin{proof}
The proof uses a similar argument as the one used in \cite{huangneely-spn}. Hence is omitted for brevity. 
\end{proof}
Note that the dual function $g(\bv{\upsilon}, \bv{\nu})$ does not contain the terms $\vartheta_k,  \varrho_k^{(s_i)}, \varphi_k^{\bv{h}_i}$. This not only simplifies the evaluation of the dual function, but also enables us to analyze the performance of our algorithm using Theorem \ref{theorem:upper-bound}. 
In the following, it is also useful to define the function $g_{s_i, \bv{h}_j}(\bv{\upsilon}, \bv{\nu})$ for each $(s_i, \bv{h}_j)$ pair: 
\begin{eqnarray}
\hspace{-.3in}&&g_{s_i, \bv{h}_j}(\bv{\upsilon}, \bv{\nu}) = \sup_{r^{nc}, \bv{P}^{(s_i)}, e^{(\bv{h}_j)}_{n}}\bigg\{ V\sum_{n, c} U^{(c)}_n(r^{nc}) \label{eq:dual-function-separable}\\
\hspace{-.3in}&&\qquad \qquad    - \sum_{n}\upsilon^{(c)}_n\big[r^{nc} +\sum_{a\in\script{N}^{(in)}_n} \mu^{(c)}_{[a, n]}(s_i, \bv{P}^{(s_i)})\nonumber \\
\hspace{-.3in}&&\qquad \qquad \qquad\qquad\qquad \quad     - \sum_{b\in\script{N}^{(o)}_n}\mu^{(c)}_{[n, b]}(s_i, \bv{P}^{(s_i)})\big] \nonumber\\
\hspace{-.3in}&& \qquad \qquad \qquad \qquad\quad  - \sum_{n}\nu_n\big[ \sum_{b \in\script{N}^{(o)}_n} P_{[n, b]}^{(s_i)} -e_n^{(\bv{h}_j)} \big] \nonumber
\bigg\}.
\end{eqnarray}
That is, $g_{s_i, \bv{h}_j}(\bv{\upsilon}, \bv{\nu}) $ is the dual function of (\ref{eq:opt-obj}) when there is a single channel state $s_i$ and a single energy state $\bv{h}_j$. 
It is easy to see from (\ref{eq:dual-function}) and (\ref{eq:dual-function-separable}) that:  
\begin{eqnarray}
g(\bv{\upsilon}, \bv{\nu}) = \sum_{s_i}\pi_{s_i}\sum_{\bv{h}_j}\pi_{\bv{h}_j}g_{s_i, \bv{h}_j}(\bv{\upsilon}, \bv{\nu}). \label{eq:dual-relation}
\end{eqnarray}

\section{Engineering the queues}\label{section:alg}
In this section, we present our Energy-limited Scheduling Algorithm (ESA) for the UMEM problem. ESA is designed based on the Lyapunov optimization technique developed in \cite{huangneely-spn} and \cite{neelynowbook}. The idea of ESA is to construct a Lyapunov scheduling algorithm with \emph{perturbed} weights for determining the energy harvesting, power allocation, routing and scheduling decisions. We will show that, by carefully perturbing the weights, one can ensure that whenever we allocate power to the links, there is always enough energy in the energy queues.

To start, we first choose a \emph{perturbation} vector $\bv{\theta}=(\theta_n, n\in\script{N})$ (to be specified later). We then define a \emph{perturbed} Lyapunov function as follows:
\begin{eqnarray}
L(t)\triangleq\frac{1}{2}\sum_{n, c\in\script{N}} \big[Q^{(c)}_n(t)\big]^2+\frac{1}{2}\sum_{n\in\script{N}} \big[E_n(t)-\theta_n\big]^2.\label{eq:lyapunov-func}
\end{eqnarray}
Now denote $\bv{Z}(t)=(\bv{Q}(t), \bv{E}(t))$, and define a one-slot conditional Lyapunov drift as follows:
\begin{eqnarray}
\Delta(t) = \expect{L(t+1) - L(t)\left.|\right. \bv{Z}(t)}. 
\end{eqnarray}
Here the expectation is taken over the randomness of the channel state and the energy state,  as well as the randomness in choosing the data  admission action, the power allocation action,  the routing and scheduling action, and the energy harvesting action. We have the following lemma regarding the drift:  
\begin{lemma}\label{lemma:drift}
Under any feasible data admission action, power allocation action, routing and scheduling action, and energy harvesting action that can be implemented at time $t$, we have: 
\begin{eqnarray}
\hspace{-.3in} && \Delta(t) -V\expect{ \sum_{n, c} U^{(c)}_n(R^{(c)}_n(t)) \left.|\right.\bv{Z}(t)}\label{eq:drift1} \\
\hspace{-.3in} && \leq B - V\expect{ \sum_{n, c} U^{(c)}_n(R^{(c)}_n(t)) \left.|\right.\bv{Z}(t)} \nonumber\\
\hspace{-.3in} &&\qquad - \sum_{n}\sum_{c}Q^{(c)}_n(t) \expect{ \sum_{b\in\script{N}^{(o)}_n}\mu_{[n, b]}^{(c)}(t) \nonumber\\
\hspace{-.3in} &&\qquad\qquad\qquad\qquad\quad - \sum_{a\in\script{N}^{(in)}_n}\mu_{[a, n]}^{(c)}(t)-R^{(c)}_n(t)\left.|\right.\bv{Z}(t)}\nonumber\\
\hspace{-.3in} && \quad\, -\sum_{n\in\script{N}} (E_n(t)-\theta_n) \expect{\sum_{b\in\script{N}^{(o)}_n}P_{[n, b]}(t)  -e_n(t) \left.|\right.\bv{Z}(t)}. \nonumber
\end{eqnarray}
Here $B=N^2(\mu^2_{max}+\frac{1}{2}R_{max}^2)+\frac{N}{2}[ P_{max}^2+h_{max}^2]$. 
\end{lemma}
\begin{proof}
See Appendix A. 
\end{proof}

Now denote the left-hand side (LHS) of (\ref{eq:drift1}) as $\Delta_V(t)$, we can rearrange the terms in (\ref{eq:drift1}) to get: 
\begin{eqnarray}
\hspace{-.3in} && \Delta_V(t) \leq B +\sum_{n\in\script{N}}(E_n(t)-\theta_n)\expect{e_n(t)\left.|\right. \bv{Z}(t)}  \label{eq:drift2}\\
\hspace{-.3in} && \qquad  - \expect{ \sum_{n, c} \big[VU^{(c)}_n(R^{(c)}_n(t)) - Q^{(c)}_n(t)R^{(c)}_n(t)\big] \left.|\right. \bv{Z}(t)} \nonumber\\ 
\hspace{-.3in} && \qquad - \expect{\sum_{n} \bigg[\sum_{c}\sum_{b\in\script{N}^{(o)}_n}\mu_{[n, b]}^{(c)}(t) \big[ Q_n^{(c)}(t) - Q_b^{(c)}(t)\big] \nonumber\\
\hspace{-.3in} && \qquad \qquad \qquad \quad \quad +(E_n(t)-\theta_n)\sum_{b\in\script{N}^{(o)}_n} P_{[n, b]}(t) \bigg]\left.|\right.\bv{Z}(t)}.  \nonumber
\end{eqnarray}

We now present the ESA algorithm. The idea of the algorithm is to approximately minimize the right-hand side (RHS) of (\ref{eq:drift2}) subject to the energy-availability constraint (\ref{eq:energycond}). 
In the algorithm, we use a parameter $\gamma\triangleq R_{max}+d_{max}\mu_{max}$, which is used in the link weight definition to allow  deterministic upper bounds on queue sizes.  

\underline{\emph{Energy-limited Scheduling Algorithm (ESA):}} Initialize $\bv{\theta}$. At every slot, observe $\bv{Q}(t)$, $\bv{E}(t)$, and $\bv{S}(t)$, do: 
\begin{itemize}
\item \underline{Energy Harvesting:} At time $t$, if $E_n(t)-\theta_n<0$, perform energy harvesting and store the harvested energy,  i.e., $e_n(t)=h_n(t)$. Else set $e_n(t)=0$. 
\item \underline{Data Admission:} At every time $t$, choose $R^{(c)}_n(t)$ to be the optimal solution of the following optimization problem: 
\begin{eqnarray}
\hspace{-.4in}\max: \,\, VU^{(c)}_n(r) - Q^{(c)}_n(t)r, \,\,\, s.t. \,\,0\leq r\leq R_{max}. \label{eq:esa-admit}
\end{eqnarray}

\item  \underline{Power Allocation:} At every time $t$, define the weight of the commodity $c$ data over link $[n, b]$ as:
\begin{eqnarray}
W^{(c)}_{[n, b]}(t) \triangleq \big[ Q_n^{(c)}(t) - Q_b^{(c)}(t) - \gamma\big]^+. 
\end{eqnarray}
Then define the link weight $W_{[n, b]}(t)=\max_{c}W^{(c)}_{[n, b]}(t)$, and  
choose $\bv{P}(t)\in\script{P}^{(s_i)}$ to maximize: 
\begin{eqnarray}
\hspace{-.4in} &&G(\bv{P}(t))\triangleq\sum_{n}\bigg[\sum_{b\in\script{N}^{(o)}_n}\mu_{[n, b]}(t)W_{[n, b]}(t) \label{eq:esa-power}\\
\hspace{-.4in} &&\qquad\qquad\qquad\qquad+(E_n(t)-\theta_n)\sum_{b\in\script{N}^{(o)}_n} P_{[n, b]}(t)\bigg],\nonumber
\end{eqnarray}
subject to the energy availability constraint (\ref{eq:energycond}). 
\item \underline{Routing and Scheduling:} For every node $n$, find any $c^*\in\rmax_{c}W^{(c)}_{[n, b]}(t)$. If $W^{(c^*)}_{[n, b]}(t)>0$, set: 
\begin{eqnarray}
\mu_{[n, b]}^{(c^*)}(t) = \mu_{[n, b]}(t),
\end{eqnarray}
that is, allocate the full rate over the link $[n, b]$ to any commodity that achieves the maximum positive weight over the link. Use idle-fill if needed. 
\item \underline{Queue Update:} Update $Q^{(c)}_n(t)$ and $E_n(t)$ according to the  dynamics (\ref{eq:Qdynamic}) and (\ref{eq:Edynamic}), respectively. 
\end{itemize}
Note that ESA only requires the knowledge of the \emph{instant} channel state $\bv{S}(t)$ and the queue sizes $\bv{Q}(t)$ and $\bv{E}(t)$. \emph{It does not even require any knowledge of the energy state process $\bv{h}(t)$.} This is very useful in practice when the knowledge of the energy source is difficult to obtained. ESA is also very different from previous algorithms for energy harvesting network, e.g., \cite{opt-energy-twc10} \cite{power-harvesting-kansal07}, where statistical knowledge of the energy source is often required. Also note that if all the links do not interfere with each other, then ESA can easily be implemented in a distributed manner, where each node only has to know about the queue sizes at its neighbor nodes and can decide on the power allocation locally.

\section{Performance Analysis}\label{section:analysis}
We now present the performance results of the ESA algorithm. In the following, we first present the results under i.i.d. network randomness and give its proof in the appendix. We later extend the performance results of ESA to the case when the network randomness is Markovian. 

\subsection{ESA under I.I.D. randomness}
Here we state the performance of ESA under the case when the channel state and the energy state, i.e., $\bv{S}(t)$ and $\bv{h}(t)$, are both i.i.d. 
\begin{theorem}\label{theorem:esa}
Under the ESA algorithm with $\theta_n\triangleq \delta\beta V +  P_{max}$ for all $n$, we have the following:
\begin{enumerate}
\item[(a)] The data queues and the energy queues satisfy the following for all time:
\begin{eqnarray}
\hspace{-.3in}&&0\leq Q^{(c)}_n(t)\leq \beta V+R_{max},\,\,\,\forall\,\,(n, c), \label{eq:data-queue-bound}\\
\hspace{-.3in}&& 0\leq E_n(t)\leq\theta_n+h_{max}, \,\,\, \forall\,\, n.\label{eq:energy-queue-bound}
\end{eqnarray}
 Moreover,  when a node $n$ allocates nonzero power to any of its outgoing links, $E_n(t)\geq P_{max}$. 
\item[(b)] Let $\bv{\overline{r}}=(\overline{r}^{nc}, \forall\, (n, c))$ be the time average admitted rate vector achieved by ESA, then:
\begin{eqnarray}
U_{tot}(\overline{\bv{r}})=\sum_{n, c}U^{(c)}_n(\overline{r}^{nc}) \geq U_{tot}(\bv{r}^*) - \frac{\tilde{B}}{V},
\end{eqnarray}
where $\bv{r}^*$ is an optimal solution of the UMEM problem, and $\tilde{B}= B+N\gamma d_{max}\mu_{max}=\Theta(1)$, i.e., independent of $V$.
\end{enumerate}
\end{theorem}
\begin{proof}
See Appendix B.
\end{proof}

We note the following of Theorem \ref{theorem:esa}: 
First, we will see that Part (a) is indeed a \emph{sample path} result. Hence it holds under \emph{arbitrary} $\bv{S}(t)$ and $\bv{h}(t)$ processes. Thus they also hold when $\bv{S}(t)$ and $\bv{h}(t)$ evolve according to some finite state irreducible and aperiodic Markov chain. 
Second, by taking $\epsilon=1/V$, we see from Part (a) that the average data queue size is $O(1/\epsilon)$. 
Combining this with Part (b), we see that that ESA achieves an $[O(\epsilon), O(1/\epsilon)]$ utility-backlog tradeoff for the UMEM problem. 
Third, we see from Part (a) that the energy queue size is deterministically upper bounded by some $O(1/\epsilon)$ constant. This provides an explicit characterization of the size of the energy storage device that is needed for achieving the desired utility performance. Such explicit bounds are particularly useful  for system deployments.



\subsection{ESA under Markovian randomness}\label{section:markovian}
We now extend our results to the more general setting where the channel state $\bv{S}(t)$ and the energy state $\bv{h}(t)$ both evolve according to some finite state irreducible and aperiodic Markov chains. Note that in this case $\pi_{s_i}$ and $\pi_{\bv{h}_i}$ represent the steady state probability of the events $\{\bv{S}(t)=s_i\}$ and $\{\bv{h}(t)=\bv{h}_i\}$, respectively. 
In this case, the performance results of ESA are summarized in the following theorem: 
\begin{theorem}\label{theorem:esa-markov}
Suppose that $\bv{S}(t)$ and $\bv{h}(t)$ evolve according to some finite state irreducible and aperiodic Markov chains. Then under ESA, we have: (a) the bounds (\ref{eq:data-queue-bound}) and (\ref{eq:energy-queue-bound}) still hold. (b) the average utility is within $O(1/V)$ of  $U_{tot}(\bv{r}^*)$, i.e., 
$U_{tot}(\overline{\bv{r}})=\sum_{n, c}U^{(c)}_n(\overline{r}^{nc}) \geq U_{tot}(\bv{r}^*) - O(1/V)$. 
\end{theorem}
\begin{proof}
Part (a) follows from Theorem \ref{theorem:esa}, since it is indeed a sample-path result. The proof of the utility performance is similar to that in \cite{huangneely_qlamarkovian}, and hence is omitted for brevity. 
\end{proof}

\section{Reducing the Buffer size}\label{section:reduce-buffer}
In this section, we show that  it is possible to achieve the same $O(\epsilon)$ close-to-optimal utility performance guarantee using energy storage devices with only $O([\log(1/\epsilon)]^2)$ sizes, while guaranteeing a much smaller average data queue size, i.e.,  $O([\log(1/\epsilon)]^2)$. Our algorithm is motivated by the following theorem, which is a modified version of Theorem $2$ in \cite{huangneely_dr_tac}. 
In the theorem,  we denote $\bv{y}=(\bv{\upsilon}, \bv{\nu})$. 
\begin{theorem}\label{theorem:exp-contraction}
Suppose  that $\bv{h}(t)$ and $S(t)$ both evolve according some finite state irreducible and aperiodic Markov chain, that $\bv{y}^*=(\bv{\upsilon}^*, \bv{\nu}^*)$ is finite and unique,  that $\bv{\theta}$ is chosen such that $\theta_n+\nu^*_n>0$, $\forall\,\,n$, and that for all $\bv{y}=(\bv{\upsilon}, \bv{\nu})$ with $\bv{\upsilon}\succeq\bv{0}, \bv{\nu}\in\mathbb{R}^N$, the dual function $g(\bv{y})$ satisfies:
\begin{eqnarray}
g(\bv{y}^*)\geq g(\bv{y})+L||\bv{y}^*-\bv{y}||, \label{eq:dualpolyhedralcond}
\end{eqnarray}
for some constant $L>0$ independent of $V$. Then $\bv{y}^*=\Theta(V)$, and that 
under ESA, there exists constants $D, K, c^*=\Theta(1)$, i.e., all independent of $V$, such that for any $m\in \mathbb{R}_+$, 
\begin{eqnarray}
\script{P}^{(r)}(D, Km)&\leq& c^*e^{-m},\label{eq:prob_pmr_special}
\end{eqnarray}
where $\script{P}^{(r)}(D, Km)$ is defined:
\begin{eqnarray}
\hspace{-.3in}&&\script{P}^{(r)}(D, Km)\triangleq\limsup_{t\rightarrow\infty}\frac{1}{t}\sum_{\tau=0}^{t-1}\prob\{\mathscr{E}(\tau, m)\},\label{eq:pmr_def} 
\end{eqnarray}
with $\mathscr{E}(t, m)$ being the following deviation event:
\begin{eqnarray}
\hspace{-.3in}&&\mathscr{E}(t, m)=\{\exists\, (n, c), |Q^{(c)}_n(t)-\upsilon^{(c)*}_{n}|>D+Km\}\label{eq:deviation-event}\\
\hspace{-.3in}&&\qquad\qquad\qquad \cup\,\{ \exists\,\,n,\,\, |(E_n(t)-\theta_n)-\nu_{n}^{*}|>D+Km\}. \nonumber
\end{eqnarray}
\end{theorem}
\begin{proof}
The proof is similar in spirit to  \cite{huangneely_dr_tac} and is omitted for brevity. 
\end{proof}
Note that the finiteness and uniqueness of $\bv{y}^*$ can usually be satisfied in practice, particularly when a certain ``slackness'' condition is met. Also note that the condition (\ref{eq:dualpolyhedralcond}) can typically be satisfied in practice when the action space is finite (See \cite{huangneely_dr_tac} for further discussions of these conditions). In this case, 
Theorem \ref{theorem:exp-contraction} states that the queue backlog vector pair is  ``exponentially attracted'' to the fixed point $(\bv{\upsilon}^*, \bv{\nu}^*+\bv{\theta})=\Theta(V)$, in that the probability of deviating decreases exponentially with the deviation distance. Therefore, the probability of deviating by some $\Theta(\log(V))$ distance will be $1/V$, which will be very small when $V$ is large. 
This suggests that most of the queue backlogs are kept in the queues to maintain a ``proper'' queue vector value to base the decisions on. If we can somehow learn the value of this vector, then we can ``subtract out'' a large amount of backlog from the network and reduce the required buffer sizes. Below, we present the Modified-ESA (MESA) algorithm. 

To start, for a given $\epsilon$, we let $V=1/\epsilon$, and define $M=4[\log(V)]^2$. 
We then associate with each node $n$ a \emph{virtual} energy queue process $\hat{E}_n(t)$ and a set of \emph{virtual} data queues $\hat{Q}^{(c)}_n(t)$ $\forall\,c$. We also associate with each node $n$ an \emph{actual} energy queue with size $M$. We assume that $V$ is chosen to be such that $\frac{M}{2}>\alpha_{max}\triangleq\max[P_{max}, h_{max}]$. 
MESA consists of two phases: Phase I runs the system using the virtual queue processes, to discover the ``attraction point'' values of the queues (as explained below). Phase II then uses these values to carefully perform the energy harvesting, power allocation, and routing and scheduling actions so as to ensure energy availability and reduce network delay. 


\underline{\emph{Modified-ESA (MESA):}} Initialize $\bv{\theta}$. Perform the following: 
\begin{itemize}
\item \underline{Phase I:} Choose a sufficiently large $T$. From time $t=0, ..., T$, run ESA using $\hat{\bv{Q}}(t)$ and  $\hat{\bv{E}}(t)$ as the data and energy queue processes. Obtain the two vectors $\bv{\script{Q}}=(\script{Q}^{(c)}_n, \,\forall\,(n,c))$ and $\bv{\script{E}}=(\script{E}_n,\forall\,n)$ by having:
\begin{eqnarray}
\script{Q}^{(c)}_n=[\hat{Q}^{(c)}_n(T)-\frac{M}{2}]^+, \,\,\script{E}=[\hat{E}_n(T) -\frac{M}{2}]^+. 
\end{eqnarray}

\item \underline{Phase II:} Reset $t=0$. 
Initialize $\hat{\bv{E}}(0)= \bv{\script{E}}$ and $\hat{\bv{Q}}(0)=\bv{\script{Q}}$. Also set $\bv{Q}(0)=\bv{0}$ and $\bv{E}(0)=\bv{0}$. 
In every time slot, first run the ESA algorithm based on $\hat{\bv{Q}}(t)$, $\hat{\bv{E}}(t)$, and $\bv{S}(t)$, to obtain the action variables, i.e., the corresponding $e_n(t)$, $R^{(c)}_n(t)$, and $\mu^{(c)}_{[n, b]}(t)$ values. 
Perform Data Admisson, Power Allocation, and Routing and Scheduling exactly as ESA, plus the following: 

\begin{itemize}
\item \underline{Energy harvesting:} If  $\hat{E}_n(t)< \script{E}_n$, let $\tilde{e}_n(t)=[e_n(t) - (\script{E}_n-\hat{E}_n(t))]^+$ and harvest $\tilde{e}(t)$ amount of energy, i.e., update $E_n(t)$ as follows:
\begin{eqnarray*}
\hspace{-.7in}&&E_n(t+1) = \big([E_n(t)   -\sum_{b\in\script{N}^{(o)}_n}P_{[n, b]}(t)]^++ \tilde{e}_n(t)\big) \land M.
\end{eqnarray*}
Here $a\land b=\min[a, b]$. Else if $\hat{E}_n(t)> \script{E}_n+M$, do not spend any power and update $E_n(t)$ according to: 
\begin{eqnarray*}
\hspace{-.3in}&&E_n(t+1) = \min\big[ E_n(t)  + e_n(t), M\big]. 
\end{eqnarray*}
Else update $E_n(t)$ according to:
\begin{eqnarray*}
\hspace{-.7in}&&E_n(t+1) = \big( [E_n(t)   -\sum_{b\in\script{N}^{(o)}_n}P_{[n, b]}(t)]^++ e_n(t)\big)\land M.
\end{eqnarray*}


\item  \underline{Packet Dropping:} 
For any node $n$ with $\hat{E}_n(t)< \script{E}_n+P_{max}$  or $\hat{E}_n(t)>\script{E}_n+M$, drop all the packets that should have been transmitted,  
i.e., change the input into any $Q_n^{(c)}(t)$ to:
\[ A^{(c)}_n(t)=R^{(c)}_{n}(t) +\sum_{a\in\script{N}^{(in)}_n}\mu^{(c)}_{[a, n]}(t)1_{[F_a(t)]}.\]
%
Here $1_{[\cdot]}$ is the indicator function and $F_a(t)$ is the event that $\hat{E}_a(t)\in [\script{E}_n+P_{max}, \script{E}+M]$. 
Then further modify the routing and scheduling action under ESA as follows: 
\begin{itemize}
\item If $\hat{Q}^{(c)}_n(t)<\script{Q}^{(c)}_n$, let $\tilde{A}^{(c)}_n(t)=\big[A^{(c)}_n(t)-[\script{Q}^{(c)}_n-\hat{Q}^{(c)}_n(t)]\big]^+$, 
update $Q^{(c)}_n(t)$ according to: 
\begin{eqnarray}
\hspace{-.6in}&& Q^{(c)}_n(t+1) \leq \big[Q^{(c)}_n(t) - \sum_{b\in\script{N}^{(o)}_n}\mu_{[n, b]}^{(c)}(t)\big]^+  + \tilde{A}_{n}^{(c)}(t).\nonumber
\end{eqnarray}

\item If $\hat{Q}^{(c)}_n(t)\geq \script{Q}_n^{(c)}$, update $Q^{(c)}_n(t)$ according to:
\begin{eqnarray}
\hspace{-.6in}&& Q^{(c)}_n(t+1) \leq \big[Q^{(c)}_n(t) - \sum_{b\in\script{N}^{(o)}_n}\mu_{[n, b]}^{(c)}(t)\big]^+  + A_{n}^{(c)}(t).\nonumber
\end{eqnarray}

\end{itemize}
\item Update $\hat{\bv{E}}(t)$ and $\hat{\bv{Q}}(t)$ using (\ref{eq:Edynamic}) and (\ref{eq:Qdynamic}). 

\end{itemize}
\end{itemize}

Note here we have used the $[\cdot]^+$ operator for updating $E_n(t)$ in the energy harvesting part. This is due to the fact that the power allocation decisions are now made based on $\hat{\bv{E}}(t)$ but not $\bv{E}(t)$. 
Note that if $\hat{E}_n(t)$ never gets below $\script{E}_n$ or above $\script{E}_n+M$, then we always have $E_n(t)=\hat{E}_n(t)-\script{E}_n$. Similarly, if $\hat{Q}^{(c)}_n(t)$ is always above $\script{Q}^{(c)}_n$ and $\hat{E}_n(t)$ is always in $[\script{E}_n+P_{max}, \script{E}_n+M]$, then we always have $Q^{(c)}_n(t)=\hat{Q}^{(c)}_n(t)-\script{Q}^{(c)}_n$. 
Our algorithm is designed to ensure that $\hat{Q}^{(c)}_n(t)$ and $\hat{E}_n(t)$ mostly stay in the ``right'' arrange, as shown in the following lemma. 
%


\begin{lemma}\label{eq:samplepath-q-bound}
For all time $t$, we have the following:
\begin{eqnarray}
\hspace{-.3in}&&   0 \leq Q^{(c)}_n(t)\leq [\hat{Q}_n^{(c)}(t) - \script{Q}^{(c)}_n]^++\gamma, \forall\, (n,c), \label{eq:Q-sp-bounds}\\
\hspace{-.3in}&&\min\big[[\hat{E}_{n}(t)-\script{E}_n]^+, M\big]\leq E_n(t). \label{eq:E-sp-ubounds}
\end{eqnarray}
\end{lemma}
\begin{proof}
See Appendix C. 
\end{proof}


We now summarize the performance result of MESA in the following theorem:
\begin{theorem}\label{theorem:mesa}
Suppose that the conditions in Theorem \ref{theorem:exp-contraction} hold, that the system is in steady state at time $T$, and that a steady state distribution for the queues  exists under ESA. 
Then under MESA with a sufficiently large $V$, with probability $1-O(\frac{1}{V^4})$, we have: 
\begin{eqnarray}
\overline{\bv{Q}}&\leq& O([\log(V)]^2), \label{eq:Q-mesa}\\
U_{tot}(\overline{\bv{r}}) &\geq& U_{tot}(\bv{r}^*) - O(1/V).\label{eq:U-mesa}
\end{eqnarray}
\end{theorem}
\begin{proof}
See Appendix D. 
\end{proof}
Note that the conditions in Theorem \ref{theorem:exp-contraction} are indeed the conditions needed for proving the exponential attraction result in \cite{huangneely_dr_tac}. Thus Theorem \ref{theorem:mesa} implies that if the exponential attraction result holds, which is mostly the case in practice (See \cite{huangneely_dr_tac} for more discussions), then one can significantly reduce the energy capacity needed to achieve the $O(\epsilon)$ close-to-optimal utility performance and greatly reduce the network congestion.



%


\section{Simulation}\label{section:simulation}
\begin{figure}[cht]
\centering
\includegraphics[height=1.2in, width=2.4in]{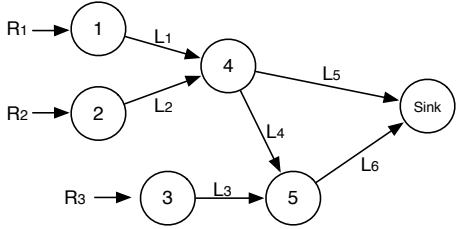}
\caption{A data collection network, where $L_i$ denotes link $i$.}\label{fig:simtopo}
\end{figure}

In this section we provide simulation results of the ESA algorithm. We consider a data collection network shown in Fig. \ref{fig:simtopo}. Such networks typically appear in the sensor network scenario where sensors are used to sense data and forward them to the sink. In this network, there are $6$ nodes. The nodes $1, 2, 3$ sense data and deliver them to node $Sink$ via the relay of nodes $4, 5$. 

\begin{figure}[cht]
\centering
\includegraphics[height=0.7in, width=1.8in]{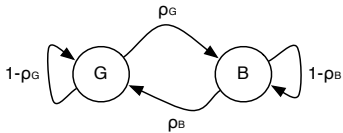}
\caption{A two-state Markov chain.}\label{fig:twostate}
\end{figure}

The channel state of each communication link $L_i$, represented by a directed edge, can be either ``G=Good'' or ``B=Bad'', and evolves according to the two-state Markov chain shown in Fig. \ref{fig:twostate} with $\rho_G=\rho_B=0.3$.  
At any time, we can allocate either zero or one unit of power. 
One unit of power can serve two packets over a link when the channel state is good, but can only serve one when the channel is bad. We assume $R_{max}=3$ and the utility functions are given by: $U_{1}(r)=U_{2}(r)=U_{3}(r)=\log(1+r)$ and $U_{4}(r)=U_{5}(r)=0$. For simplicity, we also assume that all the links do not interfere with each other. 

\begin{figure}[cht]
\centering
\includegraphics[height=1.6in, width=3.4in]{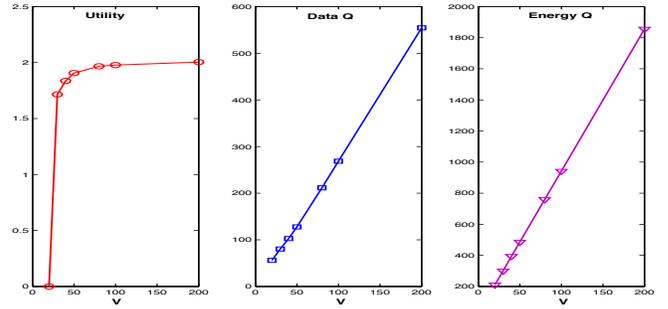}
\caption{Simulation results of ESA. }\label{fig:simplot}
\end{figure}

We also assume that for each node, the available energy $h_n(t)$ evolves according to the same two-state Markov chain in Fig. \ref{fig:twostate}. When the state is good, $h_n(t)=2$, otherwise $h_n(t)=0$. It is easy to see that in this case, $\beta=1$, $\delta=2$, $\mu_{max}=2$, $d_{max}=2$ and $P_{max}=2$. 
Using the results in Theorem \ref{theorem:esa}, we set $\theta_n=\delta\beta V+P_{max}=2V+2$. We also see that in this case, we can use $\gamma=d_{max}\mu_{max}+R_{max}=7$. The simulation results are plotted in Fig. \ref{fig:simplot}. We see in Fig. \ref{fig:simplot} that the total network utility converges quickly to very close to the optimal value, which can be shown to be roughly $2.03$. We also see that the average data queue size and the average energy queue size both grow linearly in $V$.

\begin{figure}[cht]
\centering
\includegraphics[height=2.1in, width=3.4in]{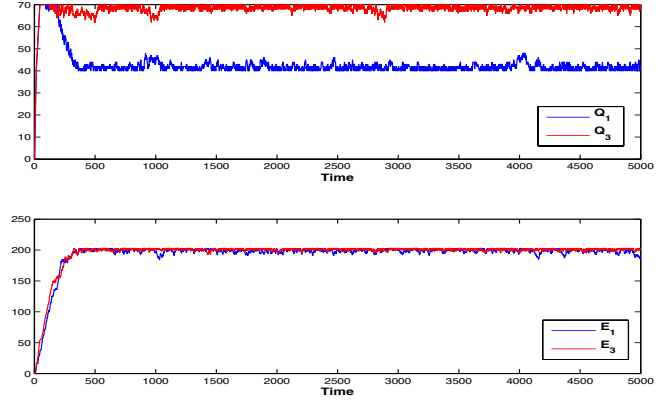}
\caption{Sample path queue processes.}\label{fig:simplot-sp}
\end{figure}

Fig. \ref{fig:simplot-sp} also shows two sample-path data queue processes and two energy queue processes under $V=100$. It can be verified that all the queue sizes satisfy the queueing bounds in Theorem \ref{theorem:esa}. Interestingly, we see that all the queue sizes are ``attracted'' to certain fixed points. 
However, different from previous work, e.g., \cite{huangneely_dr_tac}, we see that the queue size of $Q_1(t)$ does not  approach this fixed point from below. It instead first has a ``burst'' in the early time slots. This is due to the fact that the systems ``waits'' for $E_1(t)$ to come close enough to its fixed point. Such an effect can be mitigated by storing an initial energy of size $\theta$ in the energy queue. 

\begin{figure}[cht]
\centering
\includegraphics[height=1.6in, width=3.4in]{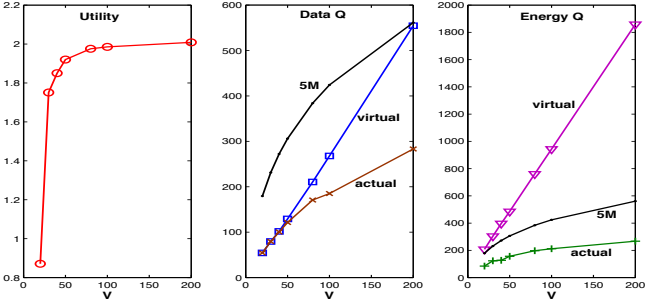}
\caption{Simulation results of MESA.}\label{fig:mesaplot}
\end{figure}

We also simulate the MESA algorithm for the same network with the same $\bv{\theta}$ value. We use $T=50V$ in Phase I for obtaining the vectors $\bv{\script{E}}$ and $\bv{\script{Q}}$. 
Fig.  \ref{fig:mesaplot} plots the performance results. We observe that extremely few packets were dropped in the simulations (at most $5$ out of more than $10^5$  packets were dropped under any $V$ values). The utility again quickly converges to the optimal as $V$ increases. We also see from the second and third plots that the actual queues only grow poly-logarithmically in $V$, i.e., $O([\log(V)]^2)$, while the virtual queues, which are the same as the actual queues under ESA, grows linearly in $V$. This shows a good match between the simulation results and Theorem \ref{theorem:mesa}. 

\section{Conclusion}\label{section:conclusion}
In this paper, we develop the Energy-limited Scheduling Algorithm (ESA) for achieving optimal utility in general energy harvesting networks equipped with only finite capacity energy storage device. We show that ESA is able to achieve an average utility that is within $O(\epsilon)$ of the optimal for any $\epsilon>0$ using energy storage devices of $O(1/\epsilon)$ sizes, while guaranteeing that the time average network congestion is $O(1/\epsilon)$. We then also develop the Modified-ESA algorithm (MESA), and show that MESA can achieve the same $O(\epsilon)$ utility performance using energy storage devices of only $O([\log(1/\epsilon)]^2)$ sizes.

\section*{Appendix A -- Proof of Lemma \ref{lemma:drift}}
Here we prove Lemma \ref{lemma:drift}
\begin{proof}
First by squaring both sides of  (\ref{eq:Qdynamic}), and using the fact that  for any $x\in\mathbb{R}$, $([x]^+)^2\leq x^2$, we have: 
\begin{eqnarray}
\hspace{-.3in} &&[Q^{(c)}_n(t+1)]^2-  [Q^{(c)}_n(t)]^2\label{eq:s1}\\
\hspace{-.3in} &&\,\,\,\leq   [\sum_{b\in\script{N}^{(o)}_n}\mu^{(c)}_{[n, b]}(t)]^2 +[\sum_{a\in\script{N}^{(in)}_n}\mu^{(c)}_{[a, n]}(t) + R^{(c)}_n(t)]^2\nonumber\\
\hspace{-.3in} &&\quad  - 2Q^{(c)}_n(t)\big[ \sum_{b\in\script{N}^{(o)}_n}\mu^{(c)}_{[n, b]}(t) - \sum_{a\in\script{N}^{(in)}_n}\mu^{(c)}_{[a, n]}(t) -R^{(c)}_n(t) \big]. \nonumber
\end{eqnarray}
By defining $\hat{B}=\frac{3}{2}d_{max}^2\mu^2_{max}+R_{max}^2$, we see that:
\begin{eqnarray}
\hspace{-.3in} &&\frac{1}{2}\big([Q^{(c)}_n(t+1)]^2 - [Q^{(c)}_n(t)]^2\big)\leq \hat{B}\label{eq:s1b}\\
\hspace{-.3in} &&\quad- Q^{(c)}_n(t)\big[ \sum_{b\in\script{N}^{(o)}_n}\mu^{(c)}_{[n, b]}(t) - \sum_{a\in\script{N}^{(in)}_n}\mu^{(c)}_{[a, n]}(t) -R^{(c)}_n(t) \big]. \nonumber
\end{eqnarray}
Using a similar approach, we get that:
\begin{eqnarray}
\hspace{-.3in} &&\frac{1}{2}\big([E_n(t+1)-\theta_n]^2 - [E_n(t)-\theta_n]^2\big)\label{eq:s2}\\
\hspace{-.3in} &&\qquad\qquad \leq \hat{B}'- [E_n(t)-\theta_n]\big[ \sum_{b\in\script{N}^{(o)}_n}P_{[n, b]}(t)  -e_n(t) \big], \nonumber
\end{eqnarray}
where $\hat{B}' = \frac{1}{2}(P_{max}+h_{max})^2$. 
Now by summing (\ref{eq:s1b}) over all $(n, c)$ and  (\ref{eq:s2}) over all $n$, and by defining $B=N^2\hat{B}+N\hat{B}'=N^2(\frac{3}{2}d_{max}^2\mu^2_{max}+R_{max}^2)+\frac{1}{2}N( P_{max}+h_{max})^2$, we have:
\begin{eqnarray*}
\hspace{-.3in} &&L(t+1)-L(t) \leq B   - \sum_{n, c}Q^{(c)}_n(t)\big[ \sum_{b\in\script{N}^{(o)}_n}\mu^{(c)}_{[n, b]}(t)\\
\hspace{-.3in} && \qquad\qquad\qquad \qquad \qquad\quad  - \sum_{a\in\script{N}^{(in)}_n}\mu^{(c)}_{[a, n]}(t) -R^{(c)}_n(t) \big]\\\
\hspace{-.3in} &&  \qquad \qquad \quad - \sum_{n}[E_n(t)-\theta_n]\big[ \sum_{b\in\script{N}^{(o)}_n}P_{[n, b]}(t)  -e_n(t) \big]. 
\end{eqnarray*}
Taking expectations on both sides over the random channel and energy states and the randomness over actions conditioning on $\bv{Z}(t)$,  subtracting from both sides the term $V\expect{\sum_{n, c} U^{(c)}_n(R^{(c)}_n(t))\left.|\right. \bv{Z}(t)}$, and rearranging the terms, we see that   the lemma follows. 
\end{proof}

\section*{Appendix B -- Proof of Theorem \ref{theorem:esa}}
Here we prove Theorem \ref{theorem:esa}. The proof idea is as follows:  
We first show that by our choice of $\bv{\theta}$, the  ESA algorithm ensures the energy-availability constraint (\ref{eq:energycond}) even if we remove it from the algorithm. This enables us to show that ESA approximately minimizes  the value of the RHS  of (\ref{eq:drift2}) over all possible policies. 
We then analyze the utility performance of ESA by relating the value of  the RHS  of (\ref{eq:drift2}) under ESA to the dual function $g_{s_i, \bv{h}_j}(\bv{\upsilon}, \bv{\nu})$. 


\begin{proof}  
(Part (a))
We first prove (\ref{eq:data-queue-bound}) using a similar argument as in \cite{neelyenergy}. It is easy to see that it holds for $t=0$, since $Q^{(c)}_n(0)=0$ for all $(n, c)$. Now assume that $Q^{(c)}_n(t)\leq \beta V+R_{max}$ for all $(n, c)$ at time $t$, we want to show that it holds for time $t+1$. 
First, if node $n$ does not receive any commodity $c$ data from other nodes, then $Q^{(c)}_n(t)\leq Q^{(c)}_n(t+1)\leq \beta V+R_{max}$. Second, if node $n$ receives endogenous commodity $c$ data from any other node $b$. Then according to the ESA algorithm, we must have:  
\[Q^{(c)}_n(t)\leq Q^{(c)}_b(t)-\gamma\leq \beta V+R_{max}-\gamma.\] However, since any node can receive at most $\gamma$ commodity $c$ packets, we have $Q^{(c)}_n(t+1) \leq \beta V+R_{max}$.
Finally, if node $n$ receives exogenous packets from outside the network, then according to (\ref{eq:esa-admit}), we must have $Q^{(c)}_n(t) \leq \beta V$. Hence $Q^{(c)}_n(t+1) \leq \beta V +R_{max}$. 

Now it is also easy to see from the energy storage part of ESA that $E_n(t)\leq\theta_n+h_{max}$, which proves (\ref{eq:energy-queue-bound}). 

We now show that if $E_n(t)\leq P_{max}$, then $G(\bv{P}(t))$ will be maximized by choosing $P_{[n, b]}(t)=0$ for all $b\in\script{N}^{(o)}_n$ at node $n$. To see this, first note that since all the actual queues are upper bounded by $\beta V+R_{max}$, we have: $W_{[n, b]}(t)\leq \beta V - d_{max}\mu_{max}$ for all $[n, b]$ and for all time. 

Now let the power allocation vector that maximizes $G(\bv{P}(t))$ at time $t$ be $\bv{P}^*$ and assume that there exists some $P^*_{[n, m]}$ that is positive. We now create a new power allocation vector $\bv{P}$ by setting only $P^*_{[n, m]}=0$ in $\bv{P}^*$. We see that $\bv{P}$ is also feasible. 
Then we have the following, in which we have written $\mu_{[n, m]}(\bv{S}(t), \bv{P}(t))$ only as a function of $\bv{P}(t)$ to simplify notation: 
\begin{eqnarray*}
\hspace{-.3in} &&\quad G(\bv{P}^*) - G(\bv{P}) \\
\hspace{-.3in} && =\sum_{n}\sum_{b\in\script{N}^{(o)}_n}\big[\mu_{[n, b]}(\bv{P}^*) - \mu_{[n, b]}(\bv{P})\big]W_{[n, b]}(t) \\
\hspace{-.3in} &&\qquad\qquad\qquad\qquad\qquad\qquad\qquad+(E_n(t)-\theta_n)P^*_{[n, m]}\nonumber\\
\hspace{-.3in} && \leq \big(\mu_{[n, m]}(\bv{P}^*) - \mu_{[n, m]}(\bv{P})\big)W_{[n, m]}(t)+(E_n(t)-\theta_n)P^*_{[n, m]}\nonumber.
\end{eqnarray*}
Here in the last step we have used (\ref{eq:rate-property2}) in Property \ref{property2} of $\mu_{[n, m]}(\cdot, \bv{P})$, which implies that $\mu_{[n, b]}(\bv{P}^*) - \mu_{[n, b]}(\bv{P})\leq0$ for all $b\neq m$. 
Now suppose $E_n(t)< P_{max}$. We see then $E_n(t)-\theta_n < -\delta\beta V$. Using Property \ref{property1} and the fact that $W_{[n, b]}(t)\leq \beta V - d_{max}\mu_{max}$, the above implies:
\begin{eqnarray*}
\hspace{-.3in} && G(\bv{P}^*) - G(\bv{P})  \leq (\beta V - d_{max}\mu_{max}) \delta P^*_{[n, m]} - \delta\beta V P^*_{[n, m]}\\
\hspace{-.3in} &&\qquad\qquad\qquad\,\,\,\,< 0. 
\end{eqnarray*}
This shows that $\bv{P}^*$ cannot have been the power vector that maximizes $G(\bv{P}(t))$ if $E_n(t)< P_{max}$. Therefore $E_n(t)\geq P_{max}$ whenever node $n$ allocates any nonzero power over any of its outgoing links. Hence all the power allocation decisions are feasible. This shows that the constraint (\ref{eq:energycond}) is indeed  redundant in ESA and 
completes the proof of Part (a). 

(Part (b)) We now prove Part (b). We first show that ESA approximately minimizes  the RHS of (\ref{eq:drift2}). To see this, note from Part (A) that ESA indeed minimizes the following function at time $t$:
\begin{eqnarray}
\hspace{-.3in} &&D(t) =  \,\, \sum_{n\in\script{N}}(E_n(t)-\theta_n) e_n(t) \label{eq:drift2a}\\
\hspace{-.3in} && \qquad   -   \sum_{n, c \in\script{N}} \big[VU^{(c)}_n(R^{(c)}_n(t)) - Q^{(c)}_n(t)R^{(c)}_n(t)\big] \nonumber\\ 
\hspace{-.3in} && \qquad   - \sum_{n\in\script{N}} \bigg[\sum_{c}\sum_{b\in\script{N}^{(o)}_n}\mu_{[n, b]}^{(c)}(t) \big[ Q_n^{(c)}(t) - Q_b^{(c)}(t) - \gamma\big] \nonumber\\
\hspace{-.3in} && \qquad \qquad \quad  \qquad \qquad \qquad +(E_n(t)-\theta_n)\sum_{b\in\script{N}^{(o)}_n} P_{[n, b]}(t)\bigg],  \nonumber
\end{eqnarray}
subject to only  the constraints: $e_n(t)\in[0, h_n(t)]$, $R_n^{(c)}(t)\in[0, R_{max}]$, $\bv{P}(t)\in\script{P}^{(s_i)}$ and (\ref{eq:ratecond}), i.e., without the energy-availability constraint (\ref{eq:energycond}). 
Now define $\tilde{D}(t)$ as follows:
\begin{eqnarray}
\hspace{-.3in} &&\tilde{D}(t) =  \sum_{n\in\script{N}}(E_n(t)-\theta_n)  e_n(t)   \label{eq:drift2b}\\
\hspace{-.3in} && \qquad  - \sum_{n, c \in\script{N}} \big[VU^{(c)}_n(R^{(c)}_n(t)) - Q^{(c)}_n(t)R^{(c)}_n(t)\big] \nonumber\\ 
\hspace{-.3in} && \qquad   - \sum_{n\in\script{N}} \bigg[\sum_{c}\sum_{b\in\script{N}^{(o)}_n}\mu_{[n, b]}^{(c)}(t) \big[ Q_n^{(c)}(t) - Q_b^{(c)}(t)\big] \nonumber\\
\hspace{-.3in} && \qquad \qquad \qquad  \quad \qquad +(E_n(t)-\theta_n)\sum_{b\in\script{N}^{(o)}_n} P_{[n, b]}(t)\bigg].  \nonumber
\end{eqnarray}
Note that $\tilde{D}(t)$ is indeed the function inside the expectation on the RHS of the drift bound (\ref{eq:drift1}). It is easy to see from the above that: 
\[D(t)=\tilde{D}(t) + \sum_n\sum_c\sum_{[n, b]\in\script{N}^{(o)}_n} \mu^{(c)}_{[n, b]}(t)\gamma.\] 
Since ESA minimizes $D(t)$, we see that:
\begin{eqnarray*}
\hspace{-.3in}&& \tilde{D}^{E}(t) + \sum_n\sum_c\sum_{b\in\script{N}^{(o)}_n}\mu^{(c)E}_{[n, b]}(t)\gamma \\
\hspace{-.3in}&&\qquad\qquad \leq \tilde{D}^{ALT}(t) + \sum_n\sum_c\sum_{b\in\script{N}^{(o)}_n}\mu^{(c)ALT}_{[n, b]}(t)\gamma,
\end{eqnarray*}
where the superscript $E$ represents the ESA algorithm, and $ALT$ represents any other alternate policy. Since 
\[0\leq \sum_n\sum_c\sum_{b\in\script{N}^{(o)}_n}\mu^{(c)}_{[n, b]}(t)\gamma\leq N^2\gamma d_{max}\mu_{max},\] 
we have:
\begin{eqnarray}
\hspace{-.3in}&& \tilde{D}^{E}(t)   \leq \tilde{D}^{ALT}(t) + N^2\gamma d_{max}\mu_{max}.\label{eq:rhs-bound}
\end{eqnarray}
That is, the value of $\tilde{D}(t)$ under ESA is no greater than its value under any other alternative policy plus a constant. 
%
%
Now using the definition of $\tilde{D}(t)$,  (\ref{eq:drift1}) can be rewritten as:
\begin{eqnarray*}
\hspace{-.3in} && \Delta(t) -V\expect{ \sum_{n, c} U^{(c)}_n(R^{(c)}_n(t)) \left.|\right.\bv{Z}(t)}  \\
\hspace{-.3in} &&\qquad\qquad\qquad\qquad\qquad \leq B + \expect{ \tilde{D}^{E}(t)\left.|\right.\bv{Z}(t)}. 
\end{eqnarray*}
Using (\ref{eq:rhs-bound}), we get:
\begin{eqnarray}
\hspace{-.3in} && \Delta(t) -V\expect{ \sum_{n, c} U^{(c)}_n(R^{(c)}_n(t)) \left.|\right.\bv{Z}(t)} \label{eq:drift2c}\\
\hspace{-.3in} && \qquad\qquad\qquad\qquad\qquad \leq \tilde{B}  + \expect{ \tilde{D}^{ALT}(t)\left.|\right.\bv{Z}(t)}, \nonumber
\end{eqnarray}
where $\tilde{B}=B+N^2\gamma d_{max}\mu_{max}$. Now consider the policy that minimizes $\tilde{D}(t)$ subject to only $e_n(t)\in[0, h_n(t)]$, $0\leq R_{n}^{(c)}(t)\leq R_{max}$, $\bv{P}\in\script{P}^{(\bv{S}(t))}$ and (\ref{eq:ratecond}), and denote the value of $\tilde{D}(t)$ under this policy by $\tilde{D}^*(t)$. It is easy to see then $\tilde{D}^*(t)$ is obtained by minimizing each term in (\ref{eq:drift2b}) over the constraints. 
Hence by comparing $\tilde{D}^*(t)$ with (\ref{eq:dual-function-separable}), we see that indeed, when $S(t)=s_i$ and $\bv{h}(t)=\bv{h}_j$, 
\[\tilde{D}^{*}(t)=-g_{s_i, \bv{h}_j}(\bv{Q}(t), \bv{\theta}-\bv{E}(t)).\]
Using this fact in  (\ref{eq:drift2c}), we have under ESA that: 
 \begin{eqnarray}
\hspace{-.3in} && \Delta(t) -V\expect{ \sum_{n, c} U^{(c)}_n(R^{(c)}_n(t)) \left.|\right.\bv{Z}(t)}\label{eq:drift3} \\
\hspace{-.3in} &&\qquad\qquad\qquad \leq \tilde{B} - \expect{ g_{s_i, \bv{h}_j}(\bv{Q}(t), \bv{\theta}-\bv{E}(t)) \left.|\right.\bv{Z}(t)}. \nonumber
\end{eqnarray}
Now using (\ref{eq:dual-relation}), i.e., $g(\bv{\upsilon}, \bv{\nu}) = \sum_{s_i}\pi_{s_i}\sum_{\bv{h}_j}\pi_{\bv{h}_j}g_{s_i, \bv{h}_j}(\bv{\upsilon}, \bv{\nu})$, the above becomes:
 \begin{eqnarray}
\hspace{-.3in} && \Delta(t) -V\expect{ \sum_{n, c} U^{(c)}_n(R^{(c)}_n(t)) \left.|\right.\bv{Z}(t)}\label{eq:drift4} \\
\hspace{-.3in} && \qquad\qquad\qquad\qquad\qquad\qquad \leq \tilde{B} -   g(\bv{Q}(t), \bv{\theta}-\bv{E}(t)). \nonumber
\end{eqnarray}
By Theorem \ref{theorem:upper-bound} and Lemma \ref{lemma:duality}, we see that:
\begin{eqnarray*}
VU_{tot}(\bv{r}^*)\leq\phi^*\leq g(\bv{\upsilon}^*, \bv{\nu}^*) \leq g(\bv{Q}(t), \bv{\theta}-\bv{E}(t)).
\end{eqnarray*}
Plug this into (\ref{eq:drift4}), we get: 
 \begin{eqnarray}
\hspace{-.3in} && \Delta(t) -V\expect{ \sum_{n, c} U^{(c)}_n(R^{(c)}_n(t)) \left.|\right.\bv{Z}(t)}  \leq \tilde{B} - VU_{tot}(\bv{r}^*). \nonumber
\end{eqnarray}
Taking expectations over $\bv{Z}(t)$ and summing the above over $t=0, ..., T-1$, we have:
\begin{eqnarray*}
\hspace{-.3in} &&\expect{L(T)-L(0)} - V\sum_{t=0}^{T-1} \expect{ \sum_{n, c} U^{(c)}_n(R^{(c)}_n(t))} \\
\hspace{-.3in} &&\qquad\qquad\qquad\qquad\qquad\qquad\qquad \leq T \tilde{B} - T VU_{tot}(\bv{r}^*). 
\end{eqnarray*}
Rearranging the terms,  using the facts that $L(t)\geq0$ and $L(0)=0$, dividing both sides by $VT$, and taking the liminf as $T\rightarrow\infty$, we get:
\begin{eqnarray*}
\hspace{-.3in} &&\liminf_{T\rightarrow\infty}\frac{1}{T}\sum_{t=0}^{T-1} \expect{ \sum_{n, c} U^{(c)}_n(R^{(c)}_n(t))} \geq U_{tot}(\bv{r}^*) -  \tilde{B}/V. 
\end{eqnarray*}
Using Jensen's inequality, we see that:
\begin{eqnarray*}
\hspace{-.3in} &&\sum_{n, c} U^{(c)}_n(\liminf_{T\rightarrow\infty}\frac{1}{T}\sum_{t=0}^{T-1} \expect{ R^{(c)}_n(t)}) \geq U_{tot}(\bv{r}^*) -  \tilde{B}/V. 
\end{eqnarray*}
This completes the proof of Part (b). 
\end{proof}

\section*{Appendix C -- Proof of Lemma \ref{eq:samplepath-q-bound}}
Here we prove Lemma  \ref{eq:samplepath-q-bound}. 
\begin{proof}
We first prove (\ref{eq:Q-sp-bounds}). We first define an intermediate process $\tilde{Q}^{(c)}_n(t)$ that evolves exactly as $Q^{(c)}_n(t)$ except that it does not discard packets when $\hat{E}_n(t)< \script{E}_n+P_{max}$  or $\hat{E}_n(t)>\script{E}_n+M$. We see then $Q^{(c)}_n(t)\leq\tilde{Q}^{(c)}_n(t)$. 
Using Lemma 3 in \cite{huangneely_dr_tac}, we see that:
$\tilde{Q}^{(c)}_n(t)\leq [\hat{Q}_n^{(c)}(t) - \script{Q}^{(c)}_n]^++\gamma$. Hence $Q^{(c)}_n(t)\leq [\hat{Q}_n^{(c)}(t) - \script{Q}^{(c)}_n]^++\gamma$ and 
(\ref{eq:Q-sp-bounds}) follows.

We now look at  (\ref{eq:E-sp-ubounds}). We see that it holds at time $0$ since $0=\hat{E}_{n}(0)-\script{E}_n= E_n(0)$. Now suppose that it  holds for $t=0, 1, ...,k$. We will show that it holds for $t=k+1$. Since if $\hat{E}_n(k+1)\leq\script{E}_n$, then (\ref{eq:E-sp-ubounds})  always holds. Below, we only consider the case when $\hat{E}_n(k+1)>\script{E}_n$, i.e., 
\begin{eqnarray}
[\hat{E}_n(k+1)-\script{E}_n]^+=\hat{E}_n(k+1)-\script{E}_n.\label{eq:foo3}
\end{eqnarray}
Also note that since all the actions are made based on $\hat{\bv{Q}}(t)$ and $\hat{\bv{E}}(t)$, by Theorem \ref{theorem:esa}, we always have $\hat{E}_n(t)\geq\sum_{b\in\script{N}^{(o)}_n}P_{[n, b]}(t)$, thus: 
\begin{eqnarray}
\hat{E}_n(t+1)=\hat{E}_n(t)-\sum_{b\in\script{N}^{(o)}_n}P_{[n, b]}(t) +e_n(t).  \label{eq:etemp}
\end{eqnarray}
We consider the following three cases: 

(I) $\hat{E}_n(k)<\script{E}_n$. Since $\hat{E}_n(k+1)>\script{E}_n$,  we must have 
$\script{E}_n - \hat{E}_n(k) \leq e_n(k)$. Then according to the harvesting rule, 
\begin{eqnarray*}
\hspace{-.3in}&&\quad E_n(k+1) \\
\hspace{-.3in}&&= \min\big[ [E_n(k)-\sum_{b\in\script{N}^{(o)}_n}P_{[n, b]}(k)]^+ \\
\hspace{-.3in}&&\qquad\qquad\qquad\qquad\qquad + e_n(k) - \script{E}_n +  \hat{E}_n(k) ,M\big]\\
\hspace{-.3in}&&\geq \min\big[ [ \hat{E}_n(k)+E_n(k)-\sum_{b\in\script{N}^{(o)}_n}P_{[n, b]}(k)]^+ \\
\hspace{-.3in}&&\qquad\qquad\qquad\qquad\qquad\qquad\qquad+ e_n(k) - \script{E}_n  ,M\big]\\
\hspace{-.3in}&&\geq \min\big[  \hat{E}_n(k)-\sum_{[n, b]}P_{[n, b]}(k)  + e_n(k) - \script{E}_n  ,M\big]\\
\hspace{-.3in}&& = \min\big[ \hat{E}_n(k+1) -\script{E}_n, M\big] \\
\hspace{-.3in}&& = \min\big[ [\hat{E}_n(k+1) -\script{E}_n]^+, M\big]. 
\end{eqnarray*}
Here the first inequality uses the property of $[\cdot]^+$, and the second inequality uses $E_n(k)\geq0$ and $\hat{E}_n(k)\geq \sum_{b\in\script{N}_n^{(o)}}P_{[n, b]}(k)$. 

(II) $\hat{E}_n(k)>\script{E}_n+M$. In this case, we see by the induction assumption that $E_n(k)=M$. Now by the update rule, we see that:
\begin{eqnarray}
E_n(k+1) = \min\big[E_n(k)+e_n(k), M\big]=M. 
\end{eqnarray}
Thus (\ref{eq:E-sp-ubounds}) still holds. 

(III) $\script{E}_n\leq\hat{E}_n(k)\leq\script{E}_n+M$. We have two cases: 

(III-A) If $\hat{E}_n(k+1)-\script{E}_n\leq M$, then using (\ref{eq:foo3}) and (\ref{eq:etemp}), we have: 
\begin{eqnarray*}
&&\min\big[[\hat{E}_n(k+1) -\script{E}_n]^+, M\big]\\
&=&  \hat{E}_n(k) - \sum_{b\in\script{N}^{(o)}_n}P_{[n, b]}(k) +e_{n}(k) - \script{E}_n\\
&\leq& \min\big[[[\hat{E}_n(k)-\script{E}_n]^+ - \sum_{b\in\script{N}^{(o)}_n}P_{[n, b]}(k)]^++e_n(k), M\big]\\
&\leq&\min\big[ [E_n(k)  - \sum_{b\in\script{N}^{(o)}_n}P_{[n, b]}(k) ]^+ +e_n(k), M\big]\\
&=& E_n(k+1). 
\end{eqnarray*}
Here the first inequality uses the property of the operator $[\cdot]^+$, and the second inequality uses the induction that $E_n(k)\geq \min\big[ [\hat{E}_n(k)-\script{E}_n]^+, M\big]=[\hat{E}_n(k)-\script{E}_n]^+$. 

(III-B) If $\hat{E}_n(k+1)-\script{E}_n> M$, then we must have $\hat{E}_n(k)\geq \script{E}_n +M-\alpha_{max}$, and that $E_n(k)\geq \hat{E}_n(k)-\script{E}_n\geq M-\alpha_{max}$. Using the fact that $\frac{M}{2}\geq\alpha_{max}$, we have:
\begin{eqnarray*}
&& E_n(k+1)\\
&=& \min\big[ E_n(k) - \sum_{b\in\script{N}^{(o)}_n}P_{[n, b]}(k) +e_n(k) ,  M   \big]\\
&\geq&  \min\big[ \hat{E}_n(k) - \sum_{b\in\script{N}^{(o)}_n}P_{[n, b]}(k) +e_n(k)- \script{E}_n,  M   \big]\\
&=& \min[\hat{E}_n(k+1)-\script{E}_n, M],
\end{eqnarray*}
which implies $E_n(k+1)=M$. Thus (\ref{eq:E-sp-ubounds}) holds. 
This completes the proof of (\ref{eq:E-sp-ubounds}) and proves the lemma. 
\end{proof}

\section*{Appendix D -- Proof of Theorem \ref{theorem:mesa}}
Here we prove Theorem \ref{theorem:mesa}. 
\begin{proof}
Since a steady state distribution  for the queues exists under the ESA algorithm, we see that $\script{P}^{(r)}(D, Km)$ is the steady state probability  that event $\mathscr{E}(t, m)$ happens.  Now consider  a large $V$ value that satisfies $\frac{M}{8}=\frac{1}{2}[\log(V)]^2\geq 2D$ and $\log(V)\geq 16K$.  We have: 
\begin{eqnarray*}
\frac{\frac{1}{2}[\log(V)]^2-D}{K}\geq \frac{\frac{1}{4}[\log(V)]^2}{K} \geq 4\log(V).
\end{eqnarray*}
By using (\ref{eq:prob_pmr_special}) and the above, we see that 
\begin{eqnarray*}
\prob(\mathscr{E}(T, \frac{\frac{1}{2}[\log(V)]^2-D}{K}))\leq c^*e^{-4\log(V)}=O(1/V^4).
\end{eqnarray*}
Using the definition of $\mathscr{E}(t, m)$, we see that when $V$ is large enough, with probability $1-O(1/V^4)$, the vectors $\hat{\bv{E}}(T)$ and $\hat{\bv{Q}}(T)$ satisfy the following for all $n, c$: 
\begin{eqnarray}
|\hat{Q}^{(c)}_n(T)-\upsilon^{(c)*}_n|\leq \frac{M}{8}, 
|\hat{E}_n(T)-(\theta_n + \nu^*_n)|\leq \frac{M}{8}. \label{eq:dev-prob}
\end{eqnarray}
Using the fact that $\script{Q}^{(c)}_n=[\hat{Q}^{(c)}_n(T)-\frac{M}{2}]^+$ and $\script{E}_n=[\hat{E}_n(T)-\frac{M}{2}]^+$, (\ref{eq:dev-prob}) and the facts that $M=4[\log(V)]^2$ and $\bv{y}^*=(\bv{\upsilon}^*, \bv{\nu}^*)=\Theta(V)$ imply that, when $V$ is large enough,  with probability $1-O(1/V^4)$, we have:
\begin{eqnarray}
\hspace{-.4in}&&-\frac{3M}{8}\geq\script{Q}^{(c)}_n-\upsilon^{(c)*}_n\geq -\frac{5M}{8}, \,\forall\, (n,c)\,\, \text{s.t.}\,\, \upsilon^{(c)*}_n\neq0, \label{eq:placeholder1}\\
\hspace{-.4in}&& \script{Q}^{(c)}_n = \upsilon^{(c)*}_n,\,\forall\, (n,c)\,\, \text{s.t.}\,\, \upsilon^{(c)*}_n=0, \label{eq:placeholder2}\\
\hspace{-.4in}&&-\frac{3M}{8}\geq \script{E}_n - (\theta_n+\nu^*_n)\geq - \frac{5M}{8}, \,\forall\, n.  \label{eq:placeholder3}
\end{eqnarray}
Having established (\ref{eq:placeholder1})-(\ref{eq:placeholder3}), (\ref{eq:Q-mesa}) can now be proven using (\ref{eq:Q-sp-bounds}) in Lemma \ref{eq:samplepath-q-bound} and a same argument as in the proof of Theorem $4$  in \cite{huangneely_dr_tac}. 

Now we consider (\ref{eq:U-mesa}). By Lemma \ref{eq:samplepath-q-bound}, when $\hat{E}_n(t)\in[\script{E}_n+P_{max}, \script{E}_n+M]$, we have $E_n(t)\geq[\hat{E}_n(t)-\script{E}_n]^+\geq P_{max}$. Thus all the power allocations are valid under MESA. 
Now since at every time $t$, MESA performs ESA's data admission, and routing and scheduling actions, if there was no packet dropping, then MESA will achieve the same utility performance as  ESA. However, since all the utility functions have  bounded derivatives, to prove the utility performance of MESA, it suffices to show that  the average rate of the packets dropped is $O(\epsilon)=O(1/V)$. 

To prove this, we first see that packet dropping happens at time $t$ only when the following event  $\hat{\mathscr{E}}(t)$ happens, i.e., 
\begin{eqnarray}
\hspace{-.4in}&&\hat{\mathscr{E}}(t) = \{\exists\,\, n, \hat{E}_n(t) <\script{E}_n+P_{max}\}\label{eq:event-mesa}\\
\hspace{-.4in}&&\qquad\qquad\qquad\quad\cup\{\exists\,\, n,  \hat{E}_n(t) >\script{E}_n+M \}\nonumber\\
\hspace{-.4in}&&\qquad\qquad\qquad\qquad\qquad\quad \cup \{\exists\,\, (n, c), \hat{Q}^{(c)}_n(t)<\script{Q}^{(c)}_n\}.\nonumber
\end{eqnarray}
However, assuming (\ref{eq:placeholder1})-(\ref{eq:placeholder3}) hold,  the following event must happen for  $\hat{\mathscr{E}}(t)$ to happen:
\begin{eqnarray*}
\tilde{\mathscr{E}}(t)&=&\{\exists\,n, \, \hat{E}_n(t)< (\theta_n+\nu^*_n)-\frac{3M}{8}+P_{max}\}\\
&&\cup \{\exists\,n, \, \hat{E}_n(t) > (\theta_n+\nu^*_n)+\frac{3M}{8}\}\\
&&\cup \{\exists\,(n, c), \, \hat{Q}^{(c)}_n(t)<\upsilon^{(c)*}_n-\frac{3M}{8}\}. 
\end{eqnarray*}
Therefore $\hat{\mathscr{E}}(t)\subset\tilde{\mathscr{E}}(t)$. 
However, it is easy to see from  (\ref{eq:deviation-event}) 
that  $\tilde{\mathscr{E}}(t)\subset\mathscr{E}(t, \tilde{m})$ with $\tilde{m}=(\frac{3M}{8}-P_{max}-D)/K=(\frac{3}{2}[\log(V)]^2-P_{max}-D)/K$. 
Therefore $\hat{\mathscr{E}}(t)\subset\mathscr{E}(t, \tilde{m})$. 
Using (\ref{eq:prob_pmr_special}) again, we see that: 
\begin{eqnarray*}
\limsup_{t\rightarrow\infty} \frac{1}{t}\sum_{\tau=0}^{t-1} \prob(\hat{\mathscr{E}}(\tau)) &\leq&  \limsup_{t\rightarrow\infty} \frac{1}{t}\sum_{\tau=0}^{t-1} \prob(\mathscr{E}(\tau,  \tilde{m}))\\
&\leq& c^*e^{-(\frac{3[\log(V)]^2}{2}-P_{max}-D)/K}.
\end{eqnarray*}
Using the facts that $\frac{1}{2}[\log(V)]^2\geq 2D$ and $\log(V)\geq 16K$, we see that:
\begin{eqnarray*}
\frac{\frac{3[\log(V)]^2}{2}-D}{K}\geq \frac{\frac{5[\log(V)]^2}{4}}{K}\geq 20\log(V). 
\end{eqnarray*}
Thus we conclude that:
\begin{eqnarray*}
\limsup_{t\rightarrow\infty} \sum_{\tau=0}^{t-1} \prob(\hat{\mathscr{E}}(\tau)) &\leq& \frac{c^*e^{P_{max}/K}}{V^{20}} = O(1/V^{20}). 
\end{eqnarray*}
Since at every time slot, the total amount of packets dropped is no more than $N(\mu_{max}+R_{max})$, we see that the average rate of packets dropped is $O(1/V)$. 
This completes the proof of Theorem  \ref{theorem:mesa}. 
\end{proof}

$\vspace{-.3in}$
\bibliographystyle{unsrt}
\bibliography{../mybib}

\end{document}